\documentclass[11pt]{article}
\usepackage{amsmath,amsfonts,amsthm,amssymb,indentfirst}

\usepackage{amsmath}
\usepackage{amsthm}
\usepackage{amsfonts}
\usepackage{amssymb}
\usepackage{color}

\theoremstyle{plain}

\theoremstyle{definition}

\newtheorem*{thmcite}{Theorem}

\newtheorem*{defncite}{Definition}

\newcommand{\xym}{\ensuremath \xymatrix@1}

\newcommand{\J}{\operatorname{\ensuremath{\it J}}}

\newcommand{\U}{\operatorname{\ensuremath {\it U}}}

\begin{document}  
\bibliographystyle{alpha}

\title{A note on completeness in the theory of strongly clean rings}
\author{Alexander J. Diesl \and Thomas J. Dorsey}
\maketitle

\begin{abstract}
Many authors have investigated the behavior of strong cleanness under certain ring extensions.  In this note, we prove that if $R$ is a ring which is complete with respect to an ideal $I$ and if $x$ is an element of $R$ whose image in $R/I$ is strongly $\pi$-regular, then $x$ is strongly clean in $R$.
\end{abstract}

In this note, all rings are associative, with unity 1, which is preserved by homomorphisms.  
Recall that, following \cite{nich:str}, an element $x$
of a ring $R$ is said to be strongly clean if there is an
idempotent $e \in R$, which commutes with $x$, such that $x - e$ is a unit in $R$. If $e$ is such an idempotent, we say that 
$x$ is $e$-strongly clean.    
The element $x\in R$ is called strongly $\pi$-regular if there is a positive integer $n$ such that $x^{n+1}\in x^nR\cap Rx^n$.  Equivalently (see \cite[Section 10]{gautamdiss}), the element $x$ is strongly $\pi$-regular if and only if $x$ satisfies Fitting's Lemma as an endomorphism of $R_R$.  This implies that $x$ is strongly $\pi$-regular if and only if there is an idempotent $e\in R$ that commutes with $x$ such that $x-e$ is a unit and $exe$ is nilpotent.  We therefore see that every strongly $\pi$-regular element of a ring is strongly clean (a well-known fact that was originally shown, using different techniques, in \cite{burgessmenal}).  A ring is said to be strongly clean (respectively, strongly $\pi$-regular) if each of its elements is strongly clean (respectively, strongly $\pi$-regular).  As an aside, following a remarkable result of Dischinger (see \cite{disch:sur} or \cite[Exercise 23.5]{lam:exs}), a ring $R$ is strongly $\pi$-regular if and only if for every $x\in R$ there is a positive integer $n$ such that $x^{n+1}\in x^nR$.

There are many examples of, and results about, strongly clean rings which happen to be complete with respect to an ideal $I$ (e.g. \cite[Theorem 2.4]{cyz:tri}, \cite[Theorem 9]{cyz:2x2}, \cite[Theorem 2.7]{yz:families}, \cite[Theorem 2.10]{fy:strclean}, and \cite[Theorem 25, Corollary 26]{mnpaper}).  The aforementioned results
demonstrate that it is frequently true (though not always, as as seen in \cite[Example 45]{tnpaper}) that, if
$R/I$ is strongly clean and $R$ is $I$-adically complete, then $R$
will be strongly clean.  Moreover, the use of completeness in this
context often simplifies proofs greatly (e.g. see \cite[Theorem
  25, Corollary 26]{mnpaper} and the surrounding discussion).  
The present note continues this theme.  

Our investigation is motivated by \cite[Theorem 2.1]{powerseries} (see also \cite[Theorem 3]{cz:nonunital} for a result in the nonunital case) which states that for a ring $R$
and a ring endomorphism $\sigma$ of $R$, an element of the
skew power series ring $R[[x;\sigma]]$ is strongly clean provided that its
constant term is strongly $\pi$-regular in $R$.  The proof is rather long and technical, and makes use of the added
structure present in the power series ring.    
The theorem in this note generalizes this result with a rather simple proof whose main trick is to consider the appropriate Peirce
decomposition when refining the relevant idempotent.  

We will need one definition which is a variant of the notion of 
  ``strong $\pi$-rad cleanness'' defined in
  \cite{diesldiss}.    
  
\begin{defncite}  Let $R$ be a ring, let $I$ be an ideal of $R$, and
  let $e \in R$ be an idempotent.  We say that $x \in R$ is
  $e$-strongly $\pi$-clean (of degree $n$), with respect to $I$, if $x$ is
  $e$-strongly clean and if, in addition, the image of $exe$ is nilpotent
  (of degree $n$) in $R/I$.
We say that $x \in R$ is strongly $\pi$-clean, with respect to $I$, if there exists an
  idempotent $e\in R$ for which $x$ is $e$-strongly $\pi$-clean, with respect to $I$.  
\end{defncite}  

Clearly, if $x$ is $e$-strongly $\pi$-clean, with respect to any ideal, then $x$ is clearly $e$-strongly clean.  
The definition above, in the case when $I$ is the Jacobson radical of $R$, agrees with the 
notion of strong $\pi$-rad cleanness defined in \cite{diesldiss}; when $I = 0$, it agrees with the 
usual notion of strong $\pi$-regularity.  In what follows, we will therefore call an element strongly $\pi$-regular 
of degree $n$ if it is strongly $\pi$-clean of degree $n$ with respect to the zero ideal.  
Also, if $x$ is $e$-strongly $\pi$-clean with respect to an ideal $I$, then it is $e$-strongly $\pi$-clean 
with respect to any ideal which contains $I$.  In particular, if $x \in R$ is strongly $\pi$-regular, it is clearly 
strongly $\pi$-clean with respect to any ideal $I$.  Note, however, that requiring that $x$ is strongly $\pi$-clean with respect to $I$ is stronger than requiring that the image of $x \in R$ is strongly $\pi$-regular in $R/I$.  
The 
difference is that the idempotent in $R/I$ which witnesses the strong $\pi$-regularity of $x$ in $R/I$ need not lift 
to an idempotent of $R$, and even if it does, it need not lift to one which commutes with $x$ in $R$.
Indeed, that is precisely what is at issue in the proof of the following result.  

\begin{thmcite}
Let $R$ be a ring, complete with respect to an ideal $I$, which is necessarily contained in the Jacobson radical $J(R)$.  Let $x\in R$.  
If the image of $x$
is strongly $\pi$-regular of degree $n$ in $R/I$, 
then $x$ is strongly $\pi$-clean of degree $n$, with respect to
$I$.  In particular, $x$ is strongly clean in
$R$.    
\end{thmcite}
\begin{proof}
For $i \ge 1$, let $\pi_i\colon R \to R/I^i$ denote the natural quotient map.  These maps will simplify the exposition, since we will be dealing simultaneously with several different quotients of $R$.
We will produce a sequence of idempotents $(e_i)_{i \ge 1}$ of $R$
with the following properties:  
\begin{enumerate}  
\item[{\rm (1)}] For each $i \ge 1$, $\pi_i(x)$ is $\pi_i(e_i)$-strongly $\pi$-clean of degree
$n$, with respect to $\pi_i(I)$,
\item[{\rm (2)}] $e_i - e_{i-1} \in I^{i-1}$ for each $i\ge 2$.
\end{enumerate}
Assuming such a sequence has been constructed, then, by (2), 
$$e = \lim_{i \rightarrow \infty} e_i$$ 
specifies a well-defined element of $R$, which is idempotent 
since $$e^2 - e = \lim_{i \rightarrow \infty} (e_i^2 - e_i) = 0.$$ 
Finally, we show that $x$ is $e$-strongly $\pi$-clean, with respect to $I$, in
$R$.  Note first that $$xe - ex = \lim_{i \rightarrow \infty} (e_i x - x e_i) = 0.$$   
Since $\pi_1(x-e) = \pi_1(x) - \pi_1(e_1)$ is a
unit in $R/I$, and $I \subseteq \J(R)$, $x-e$ is a unit in $R$.  
Finally, $(e x e)^n \in I $, since $\pi_1(exe)^n=\pi_1(e_1xe_1)^n=0$ in $R/I$ by (1).

To begin, since $\pi_1(x)$ is strongly
$\pi$-regular of degree $n$ in $R/I$, there is an idempotent $z$ in
$R/I$ such that $\pi_1(x) - z \in
\U(R/I)$, $\pi_1(x) z = z\pi_1(x)$ 
and $(z \pi_1(x) z)^n = 0$.
Since $R$ is $I$-adically complete,
idempotents lift modulo $I$ (e.g., by \cite[Theorem 21.31]{lam:fc}), and we may lift $z$ to an
idempotent of $R$, which we call $e_1$.  Note that $z = \pi(e_1)$, and that condition (1) is satisfied.    

Inductively, suppose that $e_1, \ldots, e_m$ have been constructed
satisfying conditions (1) and (2) for $1 \le i \le m$.  Write $f_m=1-e_m$.
Since $\pi_m(x)$ is $\pi_m(e_m)$-strongly $\pi$-clean with respect to $I/I^m$ in $R/I^m$, 
the Peirce decomposition of $\pi_m(x_m)$ with
respect to $\pi_m(e_m)$ is  
$$\pi_m(x_m) = \begin{pmatrix} \pi_m(e_m x_m e_m) & 0 \\ 0 & \pi_m(f_m x_m f_m) \end{pmatrix},$$
where $\pi_m(f_m x_m f_m)$ is a unit in $\pi_m(f_mRf_m)$ and $\pi_m(e_mxe_m)^n\in \pi_m(I)$.  Note that, since $f_mI^mf_m\subseteq f_mJ(R)f_m=J(f_mRf_m)$ and units lift modulo the Jacobson radical, we see that $f_m x f_m$ is a unit in $f_mRf_m$.  Note further that, since 
$(e_m x e_m)^n \in I$, it follows
  that $\pi_{m+1}(e_m x e_m)$ is nilpotent in the ring $R/I^{m+1}$.  

We now consider the Peirce
decomposition of $\pi_{m+1}(x)$ with respect to the idempotent $\pi_{m+1}(e_m)$ in $R/I^{m+1}$.  The
elements $x$ and $e_m$, which commute modulo $I^m$, need not commute
modulo $I^{m+1}$, so the Peirce decomposition of $\pi_{m+1}(x)$ in $R/I^{m+1}$
need not be diagonal.  Nevertheless, if we write it as 
$$\pi_{m+1}(x) = \begin{pmatrix} a & b \\ c & d \end{pmatrix},$$
then, by the discussion in the previous paragraph, $a$ is nilpotent in $\pi_{m+1}(e_mRe_m)$, $d$ is a unit in 
$\pi_{m+1}(f_mRf_m)$, and $b,c \in \pi_{m+1}(I^m)$.       

We will now perturb the idempotent $e_m$ to an idempotent $e_{m+1}$ which
commutes with $x$ modulo $I^{m+1}$.  
Working in $R/I^{m+1}$, the Peirce decomposition of $\pi_{m+1}(e_m)$, with respect to itself, is simply
$$\begin{pmatrix} 1 & 0 \\ 0 & 0 \end{pmatrix}.$$  
Now, for any $r \in \pi_{m+1}(e_mRf_m\cap I^m)$
 and any  $s \in \pi_{m+1}(f_mRe_m\cap I^m)$, the element $\pi_{m+1}(e_m) + r + s$, whose Peirce decomposition is 
$$\begin{pmatrix} 1 & r \\ s & 0 \end{pmatrix},$$ is idempotent in
$R/I^{m+1}$ since $(I^m)^2 \subseteq I^{m+1}$.  
We seek $r$ and $s$ as above such that $$\begin{pmatrix} 1 & r \\ s &
  0 \end{pmatrix} \begin{pmatrix} a & b \\ c & d \end{pmatrix} =
\begin{pmatrix} a & b \\ c & d \end{pmatrix} \begin{pmatrix} 1 & r \\
  s & 0 \end{pmatrix}.$$ 
Since $\pi_{m+1}(I^m)^2 = 0$ in $R/I^{m+1}$ and $b,c,r,s \in \pi_{m+1}(I^m)$,
the above equation holds if and only if $ar - rd = b$ and $ds - sa =
-c$.  Let $k$ be the index of nilpotence of $a$.  Set $r = - \sum_{i=1}^{k} a^{i-1} b d^{-i}$, which belongs to
$\pi_{m+1}(e_mRf_m\cap I^m)$, and set $s = \sum_{i=1}^k d^{-i} c a^{i-1}$,
which belongs to $\pi_{m+1}(f_mRe_m\cap I^m)$.  It is easy to check (as in
\cite[Example 13]{tnpaper}), that $ds - sa = -c$ and $ar - rd = b$, and it follows
that the idempotent $\begin{pmatrix} 1 & r \\ s & 0 \end{pmatrix}$ commutes with $\pi_{m+1}(x)$
in $R/I^{m+1}$.  Lift $\begin{pmatrix} 1 & r \\ s & 0 \end{pmatrix}$
to an idempotent $e_{m+1}$ of $R$, which clearly agrees, modulo $I^m$, with $e_m$.  
Since $\pi_m(x)-\pi_m(e_m)$ is a unit in $R/I^m$ and $x - e_{m+1}$ agrees with $x -e_m$, modulo $I^m \subseteq \J(R)$, $\pi_{m+1}(x) - \pi_{m+1}(e_{m+1})$ is a unit in $R/I^{m+1}$.  
Since $e_{m+1} x e_{m+1}$ agrees with $e_1 x e_1$, modulo $I$, and $\pi_1(e_1xe_1)$ is nilpotent of index $n$ in $R/I$, $\pi_{m+1}(e_{m+1}xe_{m+1})$ is
nilpotent of index $n$ modulo $\pi_{m+1}(I)$.  We therefore conclude that $\pi_{m+1}(x)$ is
$\pi_{m+1}(e_{m+1})$-strongly $\pi$-clean of degree $n$, with respect to $\pi_{m+1}(I)$ in $R/I^{m+1}$.  
By induction, we have produced a sequence of idempotents satisfying conditions (1) and (2), which completes the proof.
\end{proof}

A few remarks are in order.  The hypothesis that the image of $x$ in $R/I$ is
strongly $\pi$-regular cannot be weakened to the hypothesis that the image of $x$ is strongly
clean.  Let $R$ be a ring with a ring endomorphism $\sigma$ such that the triangular matrix ring $\mathbb{T}_2(R)$ is strongly clean but such that $\mathbb{T}_2(R[[x; \sigma ]])$ is not strongly clean (such an example is constructed in \cite[Example 45]{tnpaper}).  View $\mathbb{T}_2(R[[x; \sigma ]]) \cong \mathbb{T}_2(R)[[x,\sigma]]$ as a skew power series ring over $\mathbb{T}_2(R)$ (using $\sigma$ to also denote the obvious extension of $\sigma$ to a ring endomorphism of $\mathbb{T}_2(R)$), and let $I$ be the ideal of $\mathbb{T}_2(R[[x; \sigma ]])$ generated by $x$.  Then $\mathbb{T}_2(R[[x; \sigma ]])/I\cong \mathbb{T}_2(R)$ is strongly clean, and $\mathbb{T}_2(R[[x; \sigma ]])$ is complete with respect to $I$, but $\mathbb{T}_2(R[[x; \sigma ]])$ is not strongly clean and is therefore certainly not strongly $\pi$-clean with respect to $I$.

There is, however, some room for
improvement on the hypotheses.  The hypothesis that the image of $x$ in $R/I$ is strongly
$\pi$-regular can be weakened to the (substantially less pleasant
hypothesis) that there be an idempotent $e\in R$ such that
$x$ is $e$-strongly clean in $R/I$ and such that for any idempotent $e'$ which agrees with $e$
modulo $I$, the maps $$l_{e'xe'} - r_{(1-e')x(1-e')}: e'(R/I^m)(1-e')
\to e'(R/I^m)(1-e')$$ and 
$$l_{(1-e')x(1-e')} - r_{e'xe'}: (1-e')(R/I^m)e' \to
(1-e')(R/I^m)e'$$ are surjective and have the additional property that the preimage of $I^{m-1}$ is
contained in $I^{m-1}$, for each $m$.  This is precisely what is used
in the proof, and is guaranteed when $(1-e)x(1-e)$ is a unit and $exe$
is nilpotent in the respective corner rings (e.g. see
\cite[Example 13]{tnpaper}).  
We have not, however, found a hypothesis which is simultaneously more
pleasant and more general, and we expect that the main use of this result will be
in the case we explicitly stated above.

\bibliography{cleanrings}

\noindent
Alexander J. Diesl\\
Department of Mathematics \\
Wellesley College \\
Wellesley, MA 02481
USA \\
Email: {\tt adiesl@wellesley.edu}\\

\noindent
Thomas J. Dorsey\\
Center for Communications Research\\
4320 Westerra Court\\
San Diego, CA 92126-1967\\
USA\\
Email: {\tt dorsey@ccrwest.org} \\

\end{document}